\documentclass[11pt]{amsart}

\usepackage{amssymb, amsmath, amsfonts, enumerate, graphs, mathrsfs}
\usepackage{hyperref}

\newcommand{\Z}{\mathbf{Z}}
\newcommand{\R}{\mathbf{R}}

\renewcommand{\phi}{\varphi}

\DeclareMathOperator*{\conv}{conv}
\newcommand\Ehr{\operatorname{Ehr}}
\renewcommand{\P}{{\mathcal P}}
\def\th{^{\text{th}}}
\newcommand{\G}{{\mathcal G}}
\newcommand{\CC}{{\mathcal C}}

\newtheorem{theorem}{Theorem}
\newtheorem{corollary}[theorem]{Corollary}
\newtheorem{proposition}[theorem]{Proposition}
\newtheorem{lemma}[theorem]{Lemma}

\newtheorem*{example}{Example}

\newtheorem*{definition}{Definition}

\title{Grid graphs, Gorenstein polytopes, and domino stackings}

\author{Matthias Beck}
\address{Department of Mathematics\\
         San Francisco State University\\
         San Francisco, CA 94132\\
         USA}
\email{beck@math.sfsu.edu}
\urladdr{http://math.sfsu.edu/beck/}

\author{Christian Haase}
\address{Fachbereich Mathematik und Informatik\\
  Freie Universit\"at Berlin\\
  Berlin, Germany} \email{christian.haase@math.fu-berlin.de}
\urladdr{http://ehrhart.math.fu-berlin.de} 

\author{Steven V Sam}
\address{Department of Mathematics\\
  Massachusetts Institute of Technology\\
  Cambridge, MA 02139\\
  USA}
\email{ssam@math.mit.edu}
\urladdr{http://math.mit.edu/~ssam/}

\begin{document}

\thanks{We thank Thomas Zaslavsky for helpful comments on an earlier
  version of this paper. Research of Matthias Beck is supported in
  part by the NSF (DMS-0810105). Research of Christian Haase is
  supported by DFG Emmy Noether fellowship HA 4383/1. Steven Sam
  thanks the research training network {\em Methods for Discrete
    Structures} and the Berkeley mathematics department for supporting
  his stay in Berlin while part of this work was done.}

\date{9 February 2009}

\keywords{magic labellings, domino tilings, perfect matching
  polytopes, rational functions, grid graphs, Gorenstein polytopes,
  Ehrhart polynomials, recurrence relations, reciprocity theorems}

\subjclass[2000]{Primary 05A15, 05C70; Secondary 52C07.}

\begin{abstract}
  We examine domino tilings of rectangular boards, which are in
  natural bijection with perfect matchings of grid graphs. This leads
  to the study of their associated perfect matching polytopes, and we
  present some of their properties, in particular, when these
  polytopes are Gorenstein. We also introduce the notion of domino
  stackings and present some results and several open questions. Our
  techniques use results from graph theory, polyhedral geometry, and
  enumerative combinatorics.
\end{abstract}

\maketitle

\section{Introduction}

Our goal is to study, from a convex geometric point of view, domino
tilings of rectangular boards. They are in natural bijection with
perfect matchings of grid graphs.  Namely, we consider their
associated perfect matching polytopes and present some of their
properties.  In particular, we characterize when these polytopes are
Gorenstein.  (We will define all these notions shortly). We also
introduce the notion of domino stackings and present some results and
several open questions.
  
The $m \times n$ {\bf grid graph} $\G(m,n)$ is defined with vertex set
$\{(i,j) \in \Z^2 : \, 0 \le i < n,\ 0 \le j < m\}$ such that two
vertices $(i,j)$ and $(i',j')$ are adjacent if and only if $|i-i'| +
|j-j'| = 1$. The $m \times n$ {\bf torus graph} $\G_T(m,n)$ consists
of the same vertex and edge set as $\G(m,n)$ with the additional edges
$\{(0,j), (n-1,j) : \, 0 \le j < m \}$ and $\{(i,0), (i,m-1) : \, 0
\le i < n \}$. We use the convention that $V(G)$ and $E(G)$ denote the
vertex set and edge set of a graph $G$.

Given a graph $G$, $M \subseteq E(G)$ is a {\bf perfect matching} if
every vertex of $G$ is incident with exactly one edge of $M$. A
generalization of a perfect matching is the notion of a {\bf magic
  labelling of sum $t$}, which is a function $E(G) \to \Z_{ \ge 0 } $
such that for each vertex $v$, the sum of the labels of the edges
incident to $v$ equals $t$. Perfect matchings are magic labellings of
sum $t=1$. A natural question is how many magic labellings of sum $t$
a given graph $G$ has. In this paper, we are interested in the case
when $G = \G(m,n)$, and denote the number of magic labellings of sum
$t$ by $T(m,n,t)$. In particular, we can fix any two parameters, let
the other vary to get a sequence, and encode this sequence in a
generating function.

The {\bf perfect matching polytope} $\P$ associated to a graph $G$ is
defined to be the convex hull in $\R^{E(G)}$ of the incidence vectors
of all perfect matchings of $G$. There is a natural identification
between points in $\P$ and weighted graphs. We denote by $\P(m,n)$ and
$\P_T(m,n)$ the perfect matching polytopes of $\G(m,n)$ and
$\G_T(m,n)$, respectively.

For an integral polytope $\P$ (i.e., the vertices of $\P$ have only
integer coordinates), let $L_{\P}(t)$ be the function that counts the
number of integer points in $t\P := \{tx : x \in \P\}$ for $t \in \Z_{
  >0 }$. It was proved by Ehrhart \cite{ehrhart} that $L_{\P}(t)$
agrees with a polynomial with constant term 1 for all positive
integers. See also \cite{ccd} for a modern treatment of Ehrhart
theory. It is easy to see (as we will show) that the magic labellings
of sum $t$ of $\G(m,n)$ correspond bijectively to the integer points
of $t\P(m,n)$. It is a well-known result (see \cite[Chapter
4.3]{stanleyec1}) that if $p(t)$ is a polynomial, then $\sum_{t \ge 0}
p(t) z^t$ evaluates to a proper rational function (i.e., the degree of
the numerator is strictly less than the degree of the denominator) of
$z$ as follows.
\[ 
\sum_{t \ge 0} T(m,n,t) z^t = \sum_{t \ge 0} L_{\P(m,n)}(t) z^t =
\frac{h(z)}{(1-z)^{d+1}} =: \Ehr_{\P(m,n)}(z) \,,
\]
for a polynomial $h(z) = h_kz^k + h_{k-1}z^{k-1} + \cdots + h_0$ where
$k \le d$ and $h_k \ne 0$. We are interested in properties of the
sequence $(h_0, h_1, \dots, h_k)$, which we call the {\bf Ehrhart
  $h$-vector} of $\P(m,n)$. For instance, Stanley showed (for general
integral polytopes) that these numbers are nonnegative integers
\cite{stanleynonnegativity}. We characterize the values $(m,n)$ for
which $\P(m,n)$ is {\bf Gorenstein}, that is, the Ehrhart $h$-vector
is palindromic ($h_j = h_{k-j}$) \cite{mirrorsymmetry, kripo}.  An
equivalent formulation for the Gorenstein property is that there
exists an integer $k$ such that $L_{\P^\circ}(t) = L_{\P}(t-k)$ for
all $t \ge k$ and $L_{\P^\circ}(t) = 0$ for $t < k$. Here
$L_{\P^\circ}(t)$ denotes the number of interior integer points of
$t\P$.  In this case, we say that $\P$ is {\bf Gorenstein of index}
$k$.

To warm up, we show in Section~\ref{basicsection} that the Gorenstein
property holds for $\P_T(m,n)$ for some values. We can conclude this
immediately from the following more general result.

\begin{proposition} \label{kregular} Let $G$ be a $k$-regular
  bipartite graph with $\#(V(G))$ even. Then the perfect matching
  polytope $\P$ of $G$ is Gorenstein of index $k$. In particular, we
  have the functional identity $L_{\P^\circ}(t) = L_{\P}(t-k)$ for $t
  \ge k$.
\end{proposition}

\begin{corollary} Assume $m \le n$ and $n>2$ is even. 
  \begin{enumerate}[{\rm (1)}]
  \item $\P_T(1,2)$ is a point.
  \item $\P_T(1,n)$ and $\P_T(2,2)$ are Gorenstein of index 2.
  \item $\P_T(2,n)$ is Gorenstein of index 3.
  \item If $m > 2$ is even, then $\P_T(m,n)$ is Gorenstein of index
    4. 
  \end{enumerate}
\end{corollary}

It is worth noting that Tagami \cite{tagami} has shown that these and
$\P_T(2,3)$ and $\P_T(2,5)$ are all of the perfect matching polytopes
of torus graphs which are Gorenstein; the latter two are Gorenstein of
index 4 and 3, respectively.

For the non-torus grid graphs, we prove the following classification
in Section~\ref{gridsection}, which is our first main result.

\begin{theorem} \label{gorensteinthm} Assume $m \le n$. The polytope
  $\P = \P(m,n)$ is Gorenstein of index $k$ if and only if one of the
  following holds:
  \begin{enumerate}[{\rm (1)}]
  \item $m = 1$ and $n$ is even, in which case $\P$ is a point,
  \item $m = 2$, in which case $k=2$ if $n=2$, and $k=3$ for $n>2$,
  \item $m = 3$ and $n$ is even, in which case $k=5$, or
  \item $m = n = 4$, in which case $k=4$.
  \end{enumerate}
\end{theorem}

Thus, in precisely these cases, $\P(m,n)$ has a palindromic Ehrhart
$h$-vector. We can say more (Section~\ref{unimodalsection}):

\begin{theorem}\label{unimodalthm} 
  If $\P(m,n)$ is Gorenstein then $\P(m,n)$ has a unimodal Ehrhart
  $h$-vector. If $m$ and $n$ are both even, then $\P_T(m,n)$ also has
  a unimodal Ehrhart $h$-vector.
\end{theorem}

(The sequence $h_0, h_1, \dots, h_k$ is {\bf unimodal} if there exists
$j$ such that $h_0 \le \cdots \le h_{j-1} \le h_j \ge h_{j+1} \ge
\cdots \ge h_k$.) Proposition \ref{dimension} in Section
\ref{basicsection} gives the dimension of $\P(m,n)$, so coupled with
this result, we have an infinite list of Gorenstein polytopes.

The second theme of this paper concerns domino tilings. A {\bf domino
  tiling} of an $m \times n$ rectangular board is a configuration of
dominos such that the entire board is covered and there are no
overlaps. An excellent survey on general tilings can be found in
\cite{tilings}. There is an obvious correspondence between domino
tilings and perfect matchings of $\G(m,n)$. A natural question is
whether there is an analogue of magic labellings of sum $t$ of
$\G(m,n)$. One possible answer is given by domino stackings. A {\bf
  domino stacking of height $t$} of an $m \times n$ rectangular board
is a collection of $t$ domino tilings piled on top of one
another. Every such domino stacking gives a magic labelling of sum $t$
of $\G(m,n)$ in a natural way. In Section~\ref{dominosection} we
establish a simple connection between magic labellings and domino
stackings by showing that this natural map is surjective via some
general results for lattice polytopes:

\begin{proposition}\label{normalprop}
  Every magic labelling of sum $t$ of $\G(m,n)$ can be realized as a
  domino stacking of height $t$ of an $m \times n$ rectangular board.
\end{proposition}

We should remark that this result is a consequence of Hall's marriage
theorem, but that it will also follow from properties of the perfect
matching polytopes.

It was shown by Klarner and Pollack \cite{dominotilings} that for
$t=1$ and fixed $m$ (or, by symmetry, fixed $n$), the generating
function $\sum_{n \ge 0} T(m,n,1) z^n$ is a proper rational function
in $z$. In Section~\ref{dominosection}, we slightly modify this proof
to obtain the same result for general $t$.

\begin{theorem} \label{recurrence} For fixed $m$ and $t$, the
  sequences $(T(m,n,t))_{n \ge 0}$ and $(T(m,n,1)^t)_{n \ge 0}$ are
  given by a linear homogeneous recurrence relation.
\end{theorem}

Finally, in Section~\ref{geometricsection}, we give a geometric
approach to domino tilings. Our main result in this section is a new,
geometric proof of the following:

\begin{theorem}[Propp \cite{dominoreciprocity}] \label{recurrthm} If
  $T(m,n,1)$ counts the number of ways to tile an $m \times n$
  rectangular board with $2 \times 1$ dominos, then the sequence
  $(T(m,n,1))_{n \ge 0}$ is given by a linear recurrence relation, and
  furthermore, this recurrence relation satisfies the reciprocity
  relation
  \[
  T(m,n,1) = (-1)^n \, T(m,-n-2,1)
  \]
  if $m \equiv 2 \bmod 4$, and
  \[
  T(m,n,1) = T(m,-n-2,1)
  \]
  otherwise.
\end{theorem}

A remark: we interpret negative arguments to a recurrence relation to
mean that the formula for the recurrence relation is run backwards. It
was shown by Propp that even if one uses nonminimal recurrence
relations, these values at negative numbers are well defined.


\section{Basic Properties of perfect matching
  polytopes} \label{basicsection}

Let $G$ be a graph. For $x \in \R^{E(G)}$ and $S \subseteq V(G)$, we
denote by $\partial(S)$ the set of edges that are incident to exactly
one vertex in $S$, and $x(\partial(S))$ is the sum of the weights of
the edges in $\partial(S)$ in the weighted graph associated with
$x$. The following theorem gives an inequality description for $\P$.

\begin{theorem}[Edmonds \cite{matchings}] \label{hyperplane} Let $G$
  be a graph with an even number of vertices. A point $x = \left( x_e
    : \, e \in E(G) \right) $ lies in the perfect matching polytope of
  $G$ if and only if
  \begin{enumerate}[{\rm (1)}]
    \item $x_e \ge 0$ for all $e \in E(G)$,
    \item $x(\partial(v)) = 1$ for all $v \in V(G)$,
    \item \label{oddsubset} $x(\partial(S)) \ge 1$ for all $S
      \subseteq V(G)$ of odd size.
    \end{enumerate}
\end{theorem}

A complete characterization of graphs for which condition
(\ref{oddsubset}) is redundant is given in \cite{perfectmatching}. In
particular, the condition is redundant for bipartite graphs.

\begin{proof}[Proof of Proposition~\ref{kregular}]
  As we just remarked, condition \eqref{oddsubset} in
  Theorem~\ref{hyperplane} is redundant for bipartite graphs. Thus
  $L_\P (t)$ counts the integer points $x = \left( x_e : \, e \in E(G)
  \right)$ that satisfy
  \begin{enumerate}[{\rm (1)}]
    \item $x_e \ge 0$ for all $e \in E(G)$,
    \item $x(\partial(v)) = t$ for all $v \in V(G)$,
  \end{enumerate}
  whereas the interior points in $\P$ are those points satisfying (1)
  with strict inequality.  In both cases, each equation in (2)
  involves exactly $k$ variables, from which we deduce that $L_{
    \P^\circ } (t) = 0$ if $1 \le t < k$ and $L_{\P^\circ}(t) =
  L_{\P}(t-k)$ for $t \ge k$.
\end{proof}

An important piece of data associated with a polytope is its
dimension. Fortunately, for grid graphs there is a simple formula for
the dimension of its perfect matching polytope. The situation for
torus graphs is complicated only in the sense that there are several
cases to consider.

\begin{proposition} \label{dimension} Suppose that $mn$ is even.
  \begin{enumerate}[{\rm (1)}]
  \item $\dim \P(m,n) = (m-1)(n-1)$.
  \item If $n > 2$ is even, then $\dim \P_T(1,n) = 1$.
  \item If $n > 2$ is even, then $\dim \P_T(2,n) = n+1$.
  \item If $n > 1$ is odd, then $\dim \P_T(2,n) = n$.
  \item If $m > 2$ and $n > 2$ are both even, then $\dim \P_T(m,n) =
    mn+1$.
  \item If $m > 2$ is even and $n > 1$ is odd, then $\dim \P_T(m,n) =
    mn$.
  \end{enumerate}
\end{proposition}

The proof of this statement follows easily from
Theorem~\ref{thm:edmonds} below. But we first need some definitions. A
graph $G$ is said to be {\bf matching covered} if every edge of $G$
belongs to a perfect matching of $G$. A {\bf brick} is a 3-connected
and bicritical graph (i.e., removing any two vertices results in a
connected graph that has a perfect matching). A method for decomposing
general graphs into bricks is given in \cite{matchingrank}, and the
number of such bricks obtained is denoted by $B(G)$. We shall not need
this machinery; for our purposes it is enough to know that $B(G) = 0$
if $G$ is bipartite, and that $B(G) = 1$ if $G$ is a brick. There is
no overlap since the definition of a brick prevents it from being
bipartite.

\begin{theorem}[Edmonds--Lov\'asz--Pulleyblank \cite{matchingrank}]
  \label{thm:edmonds}
  Let $G$ be a matching covered graph. If $\P$ is the perfect matching
  polytope of $G$, then $\dim \P = \#(E(G)) - \#(V(G)) + 1 - B(G)$.
\end{theorem}

\begin{proof}[Proof of Proposition \ref{dimension}]
  Both $\G(m,n)$ and $\G_T(m,n)$ are matching covered if $m>1$ and
  $n>1$. The number of edges of $\G(m,n)$ is $m(n-1) + n(m-1)$, and
  the number of vertices is $mn$. In the case that either $m=1$ or
  $n=1$, $\P(m,n)$ consists of a single point, which has dimension 0
  and agrees with the formula.

  When $m$ and $n$ are both even, $\G_T(m,n)$ is bipartite, so it is
  enough to count edges. For $m = 2$, the number of edges is $3n$, and
  when $m > 2$, the number of edges is $2mn$.

  In the other cases when $m$ is even and $n>1$ is odd, $\G_T(m,n)$ is
  not bipartite, but is 3-connected and bicritical if $n > 1$. So in
  this case, $\dim \P_T(m,n) = \#(E(G)) - \#(V(G))$, and we just need
  to count edges, which we've done above for $n>1$. For the case $m>2$
  and $n=1$, there are exactly two perfect matchings of $\G_T(m,1)$,
  which gives $\dim \P_T(m,1) = 1$.
\end{proof}


\section{Characterization of Grid Graphs with Gorenstein
  Polytopes} \label{gridsection}

Now we give a lemma that will be used in the proof of
Theorem~\ref{gorensteinthm}.

\begin{lemma} \label{difflemma} Let $n>2$ be even. Given a magic
  labelling of sum $t$ of $\G(3,n)$, let $a_i$, $b_i$, and $c_i$
  denote the values of the $i\th$ edge in the top, middle, and bottom
  rows, respectively. Then for $i$ even, $c_i = b_i - a_i$ and for $i$
  odd, $c_i = t - a_i + b_i$.
\end{lemma}

\begin{proof} We refer to Figure \ref{grid3xninduct}.
\begin{figure}
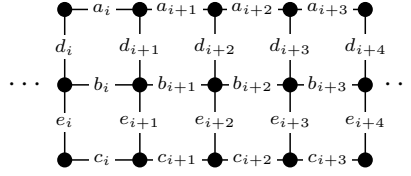

\begin{graph}(5,2)
  \freetext(0,1){$\cdots$}
  \roundnode{n1}(.5,2)
  \roundnode{n2}(1.5,2)
  \roundnode{n3}(2.5,2)
  \roundnode{n4}(3.5,2)
  \roundnode{n5}(4.5,2)
  \roundnode{n6}(.5,1)
  \roundnode{n7}(1.5,1)
  \roundnode{n8}(2.5,1)
  \roundnode{n9}(3.5,1)
  \roundnode{n10}(4.5,1)
  \roundnode{n11}(.5,0)
  \roundnode{n12}(1.5,0)
  \roundnode{n13}(2.5,0)
  \roundnode{n14}(3.5,0)
  \roundnode{n15}(4.5,0)

  \edge{n1}{n2} \edgetext{n1}{n2}{\tiny $a_i$}
  \edge{n2}{n3} \edgetext{n2}{n3}{\tiny $a_{i+1}$}
  \edge{n3}{n4} \edgetext{n3}{n4}{\tiny $a_{i+2}$}
  \edge{n4}{n5} \edgetext{n4}{n5}{\tiny $a_{i+3}$}
  \edge{n6}{n7} \edgetext{n6}{n7}{\tiny $b_i$}
  \edge{n7}{n8} \edgetext{n7}{n8}{\tiny $b_{i+1}$}
  \edge{n8}{n9} \edgetext{n8}{n9}{\tiny $b_{i+2}$}
  \edge{n9}{n10} \edgetext{n9}{n10}{\tiny $b_{i+3}$}
  \edge{n11}{n12} \edgetext{n11}{n12}{\tiny $c_i$}
  \edge{n12}{n13} \edgetext{n12}{n13}{\tiny $c_{i+1}$}
  \edge{n13}{n14} \edgetext{n13}{n14}{\tiny $c_{i+2}$}
  \edge{n14}{n15} \edgetext{n14}{n15}{\tiny $c_{i+3}$}

  \edge{n1}{n6} \edgetext{n1}{n6}{\tiny $d_i$}
  \edge{n2}{n7} \edgetext{n2}{n7}{\tiny $d_{i+1}$}
  \edge{n3}{n8} \edgetext{n3}{n8}{\tiny $d_{i+2}$}
  \edge{n4}{n9} \edgetext{n4}{n9}{\tiny $d_{i+3}$}
  \edge{n5}{n10} \edgetext{n5}{n10}{\tiny $d_{i+4}$}
  \edge{n6}{n11} \edgetext{n6}{n11}{\tiny $e_i$}
  \edge{n7}{n12} \edgetext{n7}{n12}{\tiny $e_{i+1}$}
  \edge{n8}{n13} \edgetext{n8}{n13}{\tiny $e_{i+2}$}
  \edge{n9}{n14} \edgetext{n9}{n14}{\tiny $e_{i+3}$}
  \edge{n10}{n15} \edgetext{n10}{n15}{\tiny $e_{i+4}$}

  \freetext(5,1){$\cdots$}
\end{graph}
\caption{Proof of Lemma \ref{difflemma}.}
\label{grid3xninduct}
\end{figure}
Here we are assuming that 
$i$ is even and $1 < i < n-4$.
Following the
constraints given by the vertices, we get $d_{i+2} = t - a_{i+1} -
a_{i+2}$, which implies
\begin{align*} 
  e_{i+2} &= t - (t-a_{i+1}-a_{i+2}) - b_{i+1} - b_{i+2}\\
  &= a_{i+1} + a_{i+2} - b_{i+1} - b_{i+2}\,.
\end{align*}
By induction, $c_{i+1} = t - a_{i+1} + b_{i+1}$, so
\begin{align*} 
  c_{i+2} &= t - (t-a_{i+1}+b_{i+1}) - (a_{i+1}+a_{i+2} -
  b_{i+1}-b_{i+2})\\
  &= b_{i+2} - a_{i+2}\,. 
\end{align*} 
We play the same game for $c_{i+3}$. Solving some equations gives
$d_{i+3} = t - a_{i+2} - a_{i+3}$ and
\begin{align*} 
  e_{i+3} &= t - (t-a_{i+2}-a_{i+3}) - b_{i+2} - b_{i+3}\\
  &= a_{i+2}+a_{i+3} - b_{i+2}-b_{i+3}\,. 
\end{align*} 
Using this identity, we finally conclude
\begin{align*} 
  c_{i+3} &= t - c_{i+2} - e_{i+3}\\
  &= t - (b_{i+2}-a_{i+2}) - (a_{i+2}+a_{i+3} - b_{i+2}-b_{i+3})\\
  &= t - a_{i+3} + b_{i+3}\,. 
\end{align*}
However, we have already taken care of the base case, as setting
variables with negative index to be 0 is equivalent to not having them
there at all, and the proof proceeds in the same way.
\end{proof}

Now we turn to the proof of Theorem~\ref{gorensteinthm}. Recall its
statement: $\P = \P(m,n)$ is Gorenstein if and only if one of the
following holds:
  \begin{enumerate}[{\rm (1)}]
  \item $m = 1$ and $n$ is even, 
  \item $m = 2$,
  \item $m = 3$ and $n$ is even, or 
  \item $m = n = 4$.
  \end{enumerate}

\begin{proof}[Proof of Theorem~\ref{gorensteinthm}]
  In order to show that a polytope $\P$ is Gorenstein of index $k$, we
  prove two things. First, we show that $k\P$ is the smallest dilate
  of $\P$ to contain an interior integer point, and that this point is
  unique. Second, we show that if $p$ is the unique integer point in
  $k\P^\circ$, then for any integer point $x \in t\P^\circ$, $x_i \ge
  p_i$. If our graph is bipartite, then this is enough to conclude
  that $x-p \in (t-k)\P$ using the hyperplane description given by
  Theorem~\ref{hyperplane}. This gives an obvious bijection which
  implies that $L_{\P^\circ}(t) = L_{\P}(t-k)$.

  If $m=1$ and $n$ is even, then there is exactly one matching of
  $\G(m,n)$. Then $\P$ is a point, so the Ehrhart series is
  $\Ehr_{\P}(z) = \frac{1}{1-z}$, and the statement follows.

  We now consider the case $m=2$. If $n=2$, $\P$ is a 2-regular graph
  and Proposition~\ref{kregular} applies. For $n>2$, note that the
  polytope $t\P$ has no interior points if $t<3$. For $t=3$, an
  interior point is given in Figure \ref{grid2xn}. The vertices of
  degree 3 force their adjacent edges to have weight 1, so it follows
  that this is the only such interior point. For $t \ge 3$, the edges
  of an interior point of $t\P$ must have weight $\ge 2$ if they have
  weight 2 in the figure. If not, then their adjacent edges have
  weight $t-1$ and they are adjacent to vertices of degree 3, which
  gives a contradiction. So $\P$ is Gorenstein in this case.
\begin{figure}
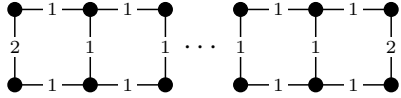

\centering 
\begin{graph}(5,1)
  \roundnode{n1}(0,1)
  \roundnode{n2}(0,0)
  \roundnode{n3}(1,1)
  \roundnode{n4}(1,0)
  \roundnode{n5}(2,1)
  \roundnode{n6}(2,0)

  \edge{n1}{n2} \edgetext{n1}{n2}{\tiny 2}
  \edge{n3}{n4} \edgetext{n3}{n4}{\tiny 1}
  \edge{n5}{n6} \edgetext{n5}{n6}{\tiny 1}
  \edge{n1}{n3} \edgetext{n1}{n3}{\tiny 1}
  \edge{n2}{n4} \edgetext{n2}{n4}{\tiny 1}
  \edge{n3}{n5} \edgetext{n3}{n5}{\tiny 1}
  \edge{n4}{n6} \edgetext{n4}{n6}{\tiny 1}
  \freetext(2.5,.5){$\cdots$}

  \roundnode{n7}(3,1)
  \roundnode{n8}(3,0)
  \roundnode{n9}(4,1)
  \roundnode{n10}(4,0)
  \roundnode{n11}(5,1)
  \roundnode{n12}(5,0)

  \edge{n7}{n8} \edgetext{n7}{n8}{\tiny 1}
  \edge{n9}{n10} \edgetext{n9}{n10}{\tiny 1}
  \edge{n11}{n12} \edgetext{n11}{n12}{\tiny 2}
  \edge{n7}{n9} \edgetext{n7}{n9}{\tiny 1}
  \edge{n8}{n10} \edgetext{n8}{n10}{\tiny 1}
  \edge{n9}{n11} \edgetext{n9}{n11}{\tiny 1}
  \edge{n10}{n12} \edgetext{n10}{n12}{\tiny 1}
\end{graph}
\caption{The interior point of $3\P$ for $\G(2,n)$.}
\label{grid2xn}
\end{figure}

For the case $m=3$, we may assume that $n>2$. The existence of
vertices of degree 4 means that $t\P$ has no interior points for
$t<4$. To get an interior point for $t=4$, consider a vertex $v$ of
degree 2. If the edges adjacent to $v$ have weight 1 and 3, then the
edge of weight 3 is adjacent to a vertex of degree 3, which cannot
happen. Thus, they must both have weight 2. However, there is a vertex
$w$ of degree 3 that is adjacent to two vertices of degree 2. Since
the edges that they share both have weight 2, there is no valid weight
for the third edge adjacent to $w$, so $4\P$ also has no interior
points. Using the notation of Lemma \ref{difflemma}, we get an
interior point of $5\P$ with the following labelling: $a_i = c_i = 3$
if $i$ is odd and $a_i = c_i = 1$ otherwise, $b_i = 1$ if $i$ is odd
and $b_i = 2$ if $i$ is even, $d_i = e_i = 2$ if $i=1$ or $i = n+1$
and $d_i = e_i = 1$ otherwise. We illustrate this interior point in
Figure \ref{grid3xninterior}.
\begin{figure}
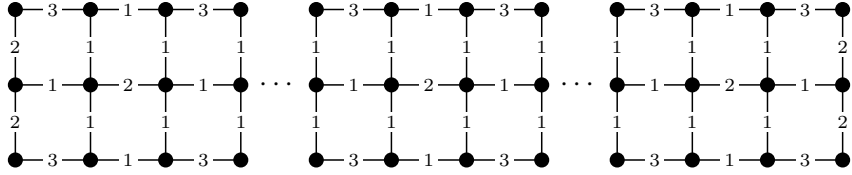

\centering 
\begin{graph}(11,2)
  \roundnode{n1}(0,2)
  \roundnode{n2}(1,2)
  \roundnode{n3}(2,2)
  \roundnode{n4}(3,2)
  \roundnode{n5}(0,1)
  \roundnode{n6}(1,1)
  \roundnode{n7}(2,1)
  \roundnode{n8}(3,1)
  \roundnode{n9}(0,0)
  \roundnode{n10}(1,0)
  \roundnode{n11}(2,0)
  \roundnode{n12}(3,0)

  \edge{n1}{n2} \edgetext{n1}{n2}{\tiny 3}
  \edge{n2}{n3} \edgetext{n2}{n3}{\tiny 1}
  \edge{n3}{n4} \edgetext{n3}{n4}{\tiny 3}
  \edge{n5}{n6} \edgetext{n5}{n6}{\tiny 1}
  \edge{n6}{n7} \edgetext{n6}{n7}{\tiny 2}
  \edge{n7}{n8} \edgetext{n7}{n8}{\tiny 1}
  \edge{n9}{n10} \edgetext{n9}{n10}{\tiny 3}
  \edge{n10}{n11} \edgetext{n10}{n11}{\tiny 1}
  \edge{n11}{n12} \edgetext{n11}{n12}{\tiny 3}

  \edge{n1}{n5} \edgetext{n1}{n5}{\tiny 2}
  \edge{n2}{n6} \edgetext{n2}{n6}{\tiny 1}
  \edge{n3}{n7} \edgetext{n3}{n7}{\tiny 1}
  \edge{n4}{n8} \edgetext{n4}{n8}{\tiny 1}
  \edge{n5}{n9} \edgetext{n5}{n9}{\tiny 2}
  \edge{n6}{n10} \edgetext{n6}{n10}{\tiny 1}
  \edge{n7}{n11} \edgetext{n7}{n11}{\tiny 1}
  \edge{n8}{n12} \edgetext{n8}{n12}{\tiny 1}

  \freetext(3.5,1){$\cdots$}

  \roundnode{n13}(4,2)
  \roundnode{n14}(5,2)
  \roundnode{n15}(6,2)
  \roundnode{n16}(7,2)
  \roundnode{n17}(4,1)
  \roundnode{n18}(5,1)
  \roundnode{n19}(6,1)
  \roundnode{n20}(7,1)
  \roundnode{n21}(4,0)
  \roundnode{n22}(5,0)
  \roundnode{n23}(6,0)
  \roundnode{n24}(7,0)

  \edge{n13}{n14} \edgetext{n13}{n14}{\tiny 3}
  \edge{n14}{n15} \edgetext{n14}{n15}{\tiny 1}
  \edge{n15}{n16} \edgetext{n15}{n16}{\tiny 3}
  \edge{n17}{n18} \edgetext{n17}{n18}{\tiny 1}
  \edge{n18}{n19} \edgetext{n18}{n19}{\tiny 2}
  \edge{n19}{n20} \edgetext{n19}{n20}{\tiny 1}
  \edge{n21}{n22} \edgetext{n21}{n22}{\tiny 3}
  \edge{n22}{n23} \edgetext{n22}{n23}{\tiny 1}
  \edge{n23}{n24} \edgetext{n23}{n24}{\tiny 3}

  \edge{n13}{n17} \edgetext{n13}{n17}{\tiny 1}
  \edge{n14}{n18} \edgetext{n14}{n18}{\tiny 1}
  \edge{n15}{n19} \edgetext{n15}{n19}{\tiny 1}
  \edge{n16}{n20} \edgetext{n16}{n20}{\tiny 1}
  \edge{n17}{n21} \edgetext{n17}{n21}{\tiny 1}
  \edge{n18}{n22} \edgetext{n18}{n22}{\tiny 1}
  \edge{n19}{n23} \edgetext{n19}{n23}{\tiny 1}
  \edge{n20}{n24} \edgetext{n20}{n24}{\tiny 1}

  \freetext(7.5,1){$\cdots$}

  \roundnode{n25}(8,2)
  \roundnode{n26}(9,2)
  \roundnode{n27}(10,2)
  \roundnode{n28}(11,2)
  \roundnode{n29}(8,1)
  \roundnode{n30}(9,1)
  \roundnode{n31}(10,1)
  \roundnode{n32}(11,1)
  \roundnode{n33}(8,0)
  \roundnode{n34}(9,0)
  \roundnode{n35}(10,0)
  \roundnode{n36}(11,0)

  \edge{n25}{n26} \edgetext{n25}{n26}{\tiny 3}
  \edge{n26}{n27} \edgetext{n26}{n27}{\tiny 1}
  \edge{n27}{n28} \edgetext{n27}{n28}{\tiny 3}
  \edge{n29}{n30} \edgetext{n29}{n30}{\tiny 1}
  \edge{n30}{n31} \edgetext{n30}{n31}{\tiny 2}
  \edge{n31}{n32} \edgetext{n31}{n32}{\tiny 1}
  \edge{n33}{n34} \edgetext{n33}{n34}{\tiny 3}
  \edge{n34}{n35} \edgetext{n34}{n35}{\tiny 1}
  \edge{n35}{n36} \edgetext{n35}{n36}{\tiny 3}

  \edge{n25}{n29} \edgetext{n25}{n29}{\tiny 1}
  \edge{n26}{n30} \edgetext{n26}{n30}{\tiny 1}
  \edge{n27}{n31} \edgetext{n27}{n31}{\tiny 1}
  \edge{n28}{n32} \edgetext{n28}{n32}{\tiny 2}
  \edge{n29}{n33} \edgetext{n29}{n33}{\tiny 1}
  \edge{n30}{n34} \edgetext{n30}{n34}{\tiny 1}
  \edge{n31}{n35} \edgetext{n31}{n35}{\tiny 1}
  \edge{n32}{n36} \edgetext{n32}{n36}{\tiny 2}
\end{graph}
\caption{The interior point of $5\P$ for $\G(3,n)$ and $n$ even.}
\label{grid3xninterior}
\end{figure}
To see this is the only interior point of $5\P$, we consider an
alternative labelling. No edge can have weight 4, so the weights for
vertices of degree 2 must be 2 and 3. Note that the 3 and 2 must be
assigned as they are in Figure \ref{grid3xninterior}. Otherwise, the
vertex of degree 3 on the far left will have sum greater than 5. We
will show that for $t \ge 5$, $a_i,c_i \ge 3$ if $i$ is odd and that
$b_i \ge 2$ if $i$ is even, which implies the uniqueness of the given
interior point. By Lemma \ref{difflemma}, we know that $c_i = b_i -
a_i$ if $i$ is even and $c_i = t - a_i + b_i$ otherwise. For the first
claim, suppose otherwise. Then $3 > t - a_i + b_i \ge t - a_i + 1$,
which implies $a_i > t - 2$, or $a_i \ge t-1$. However, $a_i$ is
adjacent to a vertex of degree 3 if $i>1$, so this is a
contradiction. For the second claim, if $b_1 = 1$, then $c_i = 1 -
a_i$, which means $c_i \le 0$, which cannot happen. Thus, we have
shown uniqueness of the interior point in $5\P$. But we have done
more. From the inequalities shown, we get for free the functional
identity $L_{\P^\circ}(t) = L_{\P}(t-k)$, giving that $\P$ is
Gorenstein.

When $m=n=4$, there are no interior points for $t\P$ when $t \le 3$,
and for $t=4$, a unique interior point is given by labelling the eight
edges incident to vertices of degree 2 with weight 2 and the rest of
the edges weight 1. For a general interior point in $t\P$ when $t \ge
4$, the edges adjacent to the corner edges must have at least weight
2, so $\P$ is Gorenstein.

For the other $(m,n)$ of interest, we split them up into two cases.
In the first case, $m \ge 4$ and $n > 4$
are both even, and in the second case, $m \ge 4$ is even and $n > 4$
is odd.

First suppose that $m \ge 4$ and $n > 4$ are both even. There are no
interior points in $t\P$ for $t < 4$. We give an interior point for
$4\P$ as follows. Assign weight 2 to the edges $\{(0, i), (0, i+1)\}$,
$\{(m-1, i), (m-1, i+1)\}$, $\{(i,0), (i+1,0)\}$, and $\{(i,n-1),
(i+1, n-1)\}$ when $i$ is even, and assign 1 to all other edges. We
illustrate this for $\G(6,6)$ in Figure
\ref{mevenneven4interior}. This example is big enough to be
instructive.
\begin{figure}
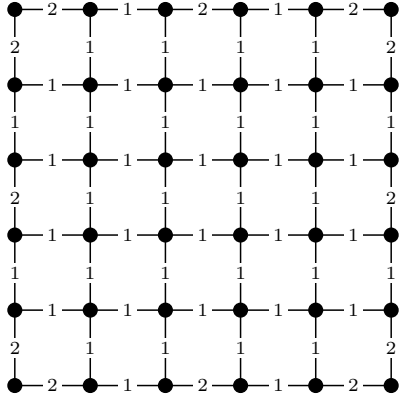

\begin{graph}(5,5)

\roundnode{n1}(0,0)
\roundnode{n2}(1,0)
\roundnode{n3}(2,0)
\roundnode{n4}(3,0)
\roundnode{n5}(4,0)
\roundnode{n6}(5,0)
\roundnode{n7}(0,1)
\roundnode{n8}(1,1)
\roundnode{n9}(2,1)
\roundnode{n10}(3,1)
\roundnode{n11}(4,1)
\roundnode{n12}(5,1)
\roundnode{n13}(0,2)
\roundnode{n14}(1,2)
\roundnode{n15}(2,2)
\roundnode{n16}(3,2)
\roundnode{n17}(4,2)
\roundnode{n18}(5,2)
\roundnode{n19}(0,3)
\roundnode{n20}(1,3)
\roundnode{n21}(2,3)
\roundnode{n22}(3,3)
\roundnode{n23}(4,3)
\roundnode{n24}(5,3)
\roundnode{n25}(0,4)
\roundnode{n26}(1,4)
\roundnode{n27}(2,4)
\roundnode{n28}(3,4)
\roundnode{n29}(4,4)
\roundnode{n30}(5,4)
\roundnode{n31}(0,5)
\roundnode{n32}(1,5)
\roundnode{n33}(2,5)
\roundnode{n34}(3,5)
\roundnode{n35}(4,5)
\roundnode{n36}(5,5)

\edge{n1}{n2} \edgetext{n1}{n2}{\tiny 2}
\edge{n2}{n3} \edgetext{n2}{n3}{\tiny 1}
\edge{n3}{n4} \edgetext{n3}{n4}{\tiny 2}
\edge{n4}{n5} \edgetext{n4}{n5}{\tiny 1}
\edge{n5}{n6} \edgetext{n5}{n6}{\tiny 2}
\edge{n7}{n8} \edgetext{n7}{n8}{\tiny 1}
\edge{n8}{n9} \edgetext{n8}{n9}{\tiny 1}
\edge{n9}{n10} \edgetext{n9}{n10}{\tiny 1}
\edge{n10}{n11} \edgetext{n10}{n11}{\tiny 1}
\edge{n11}{n12} \edgetext{n11}{n12}{\tiny 1}
\edge{n13}{n14} \edgetext{n13}{n14}{\tiny 1}
\edge{n14}{n15} \edgetext{n14}{n15}{\tiny 1}
\edge{n15}{n16} \edgetext{n15}{n16}{\tiny 1}
\edge{n16}{n17} \edgetext{n16}{n17}{\tiny 1}
\edge{n17}{n18} \edgetext{n17}{n18}{\tiny 1}
\edge{n19}{n20} \edgetext{n19}{n20}{\tiny 1}
\edge{n20}{n21} \edgetext{n20}{n21}{\tiny 1}
\edge{n21}{n22} \edgetext{n21}{n22}{\tiny 1}
\edge{n22}{n23} \edgetext{n22}{n23}{\tiny 1}
\edge{n23}{n24} \edgetext{n23}{n24}{\tiny 1}
\edge{n25}{n26} \edgetext{n25}{n26}{\tiny 1}
\edge{n26}{n27} \edgetext{n26}{n27}{\tiny 1}
\edge{n27}{n28} \edgetext{n27}{n28}{\tiny 1}
\edge{n28}{n29} \edgetext{n28}{n29}{\tiny 1}
\edge{n29}{n30} \edgetext{n29}{n30}{\tiny 1}
\edge{n31}{n32} \edgetext{n31}{n32}{\tiny 2}
\edge{n32}{n33} \edgetext{n32}{n33}{\tiny 1}
\edge{n33}{n34} \edgetext{n33}{n34}{\tiny 2}
\edge{n34}{n35} \edgetext{n34}{n35}{\tiny 1}
\edge{n35}{n36} \edgetext{n35}{n36}{\tiny 2}

\edge{n1}{n7} \edgetext{n1}{n7}{\tiny 2}
\edge{n2}{n8} \edgetext{n2}{n8}{\tiny 1}
\edge{n3}{n9} \edgetext{n3}{n9}{\tiny 1}
\edge{n4}{n10} \edgetext{n4}{n10}{\tiny 1}
\edge{n5}{n11} \edgetext{n5}{n11}{\tiny 1}
\edge{n6}{n12} \edgetext{n6}{n12}{\tiny 2}
\edge{n7}{n13} \edgetext{n7}{n13}{\tiny 1}
\edge{n8}{n14} \edgetext{n8}{n14}{\tiny 1}
\edge{n9}{n15} \edgetext{n9}{n15}{\tiny 1}
\edge{n10}{n16} \edgetext{n10}{n16}{\tiny 1}
\edge{n11}{n17} \edgetext{n11}{n17}{\tiny 1}
\edge{n12}{n18} \edgetext{n12}{n18}{\tiny 1}
\edge{n13}{n19} \edgetext{n13}{n19}{\tiny 2}
\edge{n14}{n20} \edgetext{n14}{n20}{\tiny 1}
\edge{n15}{n21} \edgetext{n15}{n21}{\tiny 1}
\edge{n16}{n22} \edgetext{n16}{n22}{\tiny 1}
\edge{n17}{n23} \edgetext{n17}{n23}{\tiny 1}
\edge{n18}{n24} \edgetext{n18}{n24}{\tiny 2}
\edge{n19}{n25} \edgetext{n19}{n25}{\tiny 1}
\edge{n20}{n26} \edgetext{n20}{n26}{\tiny 1}
\edge{n21}{n27} \edgetext{n21}{n27}{\tiny 1}
\edge{n22}{n28} \edgetext{n22}{n28}{\tiny 1}
\edge{n23}{n29} \edgetext{n23}{n29}{\tiny 1}
\edge{n24}{n30} \edgetext{n24}{n30}{\tiny 1}
\edge{n25}{n31} \edgetext{n25}{n31}{\tiny 2}
\edge{n26}{n32} \edgetext{n26}{n32}{\tiny 1}
\edge{n27}{n33} \edgetext{n27}{n33}{\tiny 1}
\edge{n28}{n34} \edgetext{n28}{n34}{\tiny 1}
\edge{n29}{n35} \edgetext{n29}{n35}{\tiny 1}
\edge{n30}{n36} \edgetext{n30}{n36}{\tiny 2}

\end{graph}
\caption{The interior point of $4\P$ for $\G(6,6)$.}
\label{mevenneven4interior}
\end{figure}
To see this is unique, note that the weights of the edges of the
corner vertices must be 2, and the weights of the edges incident to a
vertex of degree 4 must be 1. This, however, forces the rest of
the values on the edges incident to vertices of degree 3. Now consider
the following interior point of $5\P$. 
\begin{itemize}
\item Assign weight 2 to edges of the form $\{(i,0), (i,1)\}$,
  $\{(i,j), (i,j+1)\}$, and $\{(i,n-2), (i,n-1)\}$, for $i \in \{0,
  m-1\}$ and $j$ odd. 
\item Assign weight 3 to edges of the form $\{(i,j), (i+1,j)\}$ for
  $i$ even and $j \in \{0, n-1\}$. 
\item Assign weight 2 to edges of the form $\{(i,j), (i+1, j)\}$ for
  $i$ odd and $j \in \{1, n-2\}$. 
\item Assign weight 2 to edges of the form $\{(i,j), (i+1, j)\}$ for
  $i$ even and $2 \le j \le n-3$. 
\item Finally, assign weight 1 to all other edges.
\end{itemize}
We illustrate this for $\G(6,6)$ in Figure \ref{mevenneven5interior}
(here $(0,0)$ is the bottom left vertex and $(5,5)$ is the top right
vertex).
\begin{figure}
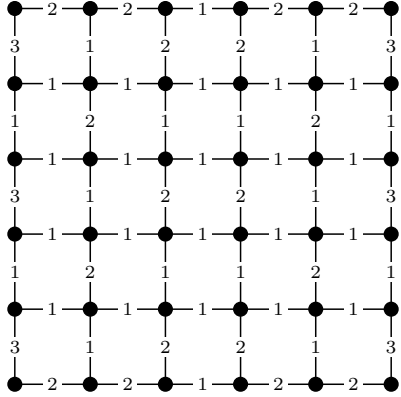

\begin{graph}(5,5)

\roundnode{n1}(0,0)
\roundnode{n2}(1,0)
\roundnode{n3}(2,0)
\roundnode{n4}(3,0)
\roundnode{n5}(4,0)
\roundnode{n6}(5,0)
\roundnode{n7}(0,1)
\roundnode{n8}(1,1)
\roundnode{n9}(2,1)
\roundnode{n10}(3,1)
\roundnode{n11}(4,1)
\roundnode{n12}(5,1)
\roundnode{n13}(0,2)
\roundnode{n14}(1,2)
\roundnode{n15}(2,2)
\roundnode{n16}(3,2)
\roundnode{n17}(4,2)
\roundnode{n18}(5,2)
\roundnode{n19}(0,3)
\roundnode{n20}(1,3)
\roundnode{n21}(2,3)
\roundnode{n22}(3,3)
\roundnode{n23}(4,3)
\roundnode{n24}(5,3)
\roundnode{n25}(0,4)
\roundnode{n26}(1,4)
\roundnode{n27}(2,4)
\roundnode{n28}(3,4)
\roundnode{n29}(4,4)
\roundnode{n30}(5,4)
\roundnode{n31}(0,5)
\roundnode{n32}(1,5)
\roundnode{n33}(2,5)
\roundnode{n34}(3,5)
\roundnode{n35}(4,5)
\roundnode{n36}(5,5)

\edge{n1}{n2} \edgetext{n1}{n2}{\tiny 2}
\edge{n2}{n3} \edgetext{n2}{n3}{\tiny 2}
\edge{n3}{n4} \edgetext{n3}{n4}{\tiny 1}
\edge{n4}{n5} \edgetext{n4}{n5}{\tiny 2}
\edge{n5}{n6} \edgetext{n5}{n6}{\tiny 2}
\edge{n7}{n8} \edgetext{n7}{n8}{\tiny 1}
\edge{n8}{n9} \edgetext{n8}{n9}{\tiny 1}
\edge{n9}{n10} \edgetext{n9}{n10}{\tiny 1}
\edge{n10}{n11} \edgetext{n10}{n11}{\tiny 1}
\edge{n11}{n12} \edgetext{n11}{n12}{\tiny 1}
\edge{n13}{n14} \edgetext{n13}{n14}{\tiny 1}
\edge{n14}{n15} \edgetext{n14}{n15}{\tiny 1}
\edge{n15}{n16} \edgetext{n15}{n16}{\tiny 1}
\edge{n16}{n17} \edgetext{n16}{n17}{\tiny 1}
\edge{n17}{n18} \edgetext{n17}{n18}{\tiny 1}
\edge{n19}{n20} \edgetext{n19}{n20}{\tiny 1}
\edge{n20}{n21} \edgetext{n20}{n21}{\tiny 1}
\edge{n21}{n22} \edgetext{n21}{n22}{\tiny 1}
\edge{n22}{n23} \edgetext{n22}{n23}{\tiny 1}
\edge{n23}{n24} \edgetext{n23}{n24}{\tiny 1}
\edge{n25}{n26} \edgetext{n25}{n26}{\tiny 1}
\edge{n26}{n27} \edgetext{n26}{n27}{\tiny 1}
\edge{n27}{n28} \edgetext{n27}{n28}{\tiny 1}
\edge{n28}{n29} \edgetext{n28}{n29}{\tiny 1}
\edge{n29}{n30} \edgetext{n29}{n30}{\tiny 1}
\edge{n31}{n32} \edgetext{n31}{n32}{\tiny 2}
\edge{n32}{n33} \edgetext{n32}{n33}{\tiny 2}
\edge{n33}{n34} \edgetext{n33}{n34}{\tiny 1}
\edge{n34}{n35} \edgetext{n34}{n35}{\tiny 2}
\edge{n35}{n36} \edgetext{n35}{n36}{\tiny 2}

\edge{n1}{n7} \edgetext{n1}{n7}{\tiny 3}
\edge{n2}{n8} \edgetext{n2}{n8}{\tiny 1}
\edge{n3}{n9} \edgetext{n3}{n9}{\tiny 2}
\edge{n4}{n10} \edgetext{n4}{n10}{\tiny 2}
\edge{n5}{n11} \edgetext{n5}{n11}{\tiny 1}
\edge{n6}{n12} \edgetext{n6}{n12}{\tiny 3}
\edge{n7}{n13} \edgetext{n7}{n13}{\tiny 1}
\edge{n8}{n14} \edgetext{n8}{n14}{\tiny 2}
\edge{n9}{n15} \edgetext{n9}{n15}{\tiny 1}
\edge{n10}{n16} \edgetext{n10}{n16}{\tiny 1}
\edge{n11}{n17} \edgetext{n11}{n17}{\tiny 2}
\edge{n12}{n18} \edgetext{n12}{n18}{\tiny 1}
\edge{n13}{n19} \edgetext{n13}{n19}{\tiny 3}
\edge{n14}{n20} \edgetext{n14}{n20}{\tiny 1}
\edge{n15}{n21} \edgetext{n15}{n21}{\tiny 2}
\edge{n16}{n22} \edgetext{n16}{n22}{\tiny 2}
\edge{n17}{n23} \edgetext{n17}{n23}{\tiny 1}
\edge{n18}{n24} \edgetext{n18}{n24}{\tiny 3}
\edge{n19}{n25} \edgetext{n19}{n25}{\tiny 1}
\edge{n20}{n26} \edgetext{n20}{n26}{\tiny 2}
\edge{n21}{n27} \edgetext{n21}{n27}{\tiny 1}
\edge{n22}{n28} \edgetext{n22}{n28}{\tiny 1}
\edge{n23}{n29} \edgetext{n23}{n29}{\tiny 2}
\edge{n24}{n30} \edgetext{n24}{n30}{\tiny 1}
\edge{n25}{n31} \edgetext{n25}{n31}{\tiny 3}
\edge{n26}{n32} \edgetext{n26}{n32}{\tiny 1}
\edge{n27}{n33} \edgetext{n27}{n33}{\tiny 2}
\edge{n28}{n34} \edgetext{n28}{n34}{\tiny 2}
\edge{n29}{n35} \edgetext{n29}{n35}{\tiny 1}
\edge{n30}{n36} \edgetext{n30}{n36}{\tiny 3}

\end{graph}
\caption{An interior point of $5\P$ for $\G(6,6)$.}
\label{mevenneven5interior}
\end{figure}
In order for $\P$ to be Gorenstein, we must have $L_{\P^\circ}(t+4) =
L_{\P}(t)$. Note that for any point of $\P$, we may add the interior
point of $4\P$ to get an interior point of $5\P$. However, all points
acquired through this must have weight $\ge 2$ for the edge $\{(0,2),
(0,3)\}$, but we have given an interior point of $5\P$ such that this
edge has weight 1, which implies $L_{\P^\circ}(5) > L_{\P}(1)$, so the
functional identity is not satisfied.

In the case that $m \ge 4$ is even and $n$ is odd, $t\P$ does not have
an interior point for $t \le 4$. To see that this is true for $t = 4$,
the edges incident to the corner vertices must have weight 2. 
This forces a labelling as in the case $n$ is even, but it will not
work because $n$ is now odd. In other words, if we label the outside
and work toward the middle, there will be no satisfactory values for
the middle edges. We now give an interior point for $5\P$.
\begin{itemize}
\item Assign weight 3 to all edges of the form $\{(i,j), (i,j+1)\}$
  where $i \in \{0,n-1\}$ and $j$ is even. 
\item Assign weight 2 to edges of the form $\{(i,j), (i+1,j)\}$ where
  $j \in \{0,m-1\}$ and $i \ne 1$. 
\item Assign weight 2 to edges $\{(i,j), (i,j+1)\}$ where $i \in
  \{1,2\}$ and $j$ is even. 
\item Assign weight 2 to edges $\{(i,j), (i,j+1)\}$ where $3 \le i \le
  n-2$ and $j$ is odd. 
\item Finally, assign weight 1 to all other edges.
\end{itemize}
This labelling is illustrated for $\G(6,7)$ in Figure
\ref{mevennodd5interior} (again, $(0,0)$ is the bottom left vertex and
$(5,6)$ is the top right).
\begin{figure}
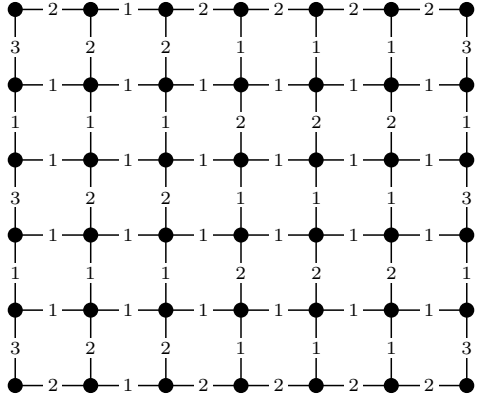

\begin{graph}(6,5)

\roundnode{n1}(0,0)
\roundnode{n2}(1,0)
\roundnode{n3}(2,0)
\roundnode{n4}(3,0)
\roundnode{n5}(4,0)
\roundnode{n6}(5,0)
\roundnode{n7}(6,0)
\roundnode{n8}(0,1)
\roundnode{n9}(1,1)
\roundnode{n10}(2,1)
\roundnode{n11}(3,1)
\roundnode{n12}(4,1)
\roundnode{n13}(5,1)
\roundnode{n14}(6,1)
\roundnode{n15}(0,2)
\roundnode{n16}(1,2)
\roundnode{n17}(2,2)
\roundnode{n18}(3,2)
\roundnode{n19}(4,2)
\roundnode{n20}(5,2)
\roundnode{n21}(6,2)
\roundnode{n22}(0,3)
\roundnode{n23}(1,3)
\roundnode{n24}(2,3)
\roundnode{n25}(3,3)
\roundnode{n26}(4,3)
\roundnode{n27}(5,3)
\roundnode{n28}(6,3)
\roundnode{n29}(0,4)
\roundnode{n30}(1,4)
\roundnode{n31}(2,4)
\roundnode{n32}(3,4)
\roundnode{n33}(4,4)
\roundnode{n34}(5,4)
\roundnode{n35}(6,4)
\roundnode{n36}(0,5)
\roundnode{n37}(1,5)
\roundnode{n38}(2,5)
\roundnode{n39}(3,5)
\roundnode{n40}(4,5)
\roundnode{n41}(5,5)
\roundnode{n42}(6,5)

\edge{n1}{n2} \edgetext{n1}{n2}{\tiny 2}
\edge{n2}{n3} \edgetext{n2}{n3}{\tiny 1}
\edge{n3}{n4} \edgetext{n3}{n4}{\tiny 2}
\edge{n4}{n5} \edgetext{n4}{n5}{\tiny 2}
\edge{n5}{n6} \edgetext{n5}{n6}{\tiny 2}
\edge{n6}{n7} \edgetext{n6}{n7}{\tiny 2}
\edge{n8}{n9} \edgetext{n8}{n9}{\tiny 1}
\edge{n9}{n10} \edgetext{n9}{n10}{\tiny 1}
\edge{n10}{n11} \edgetext{n10}{n11}{\tiny 1}
\edge{n11}{n12} \edgetext{n11}{n12}{\tiny 1}
\edge{n12}{n13} \edgetext{n12}{n13}{\tiny 1}
\edge{n13}{n14} \edgetext{n13}{n14}{\tiny 1}
\edge{n15}{n16} \edgetext{n15}{n16}{\tiny 1}
\edge{n16}{n17} \edgetext{n16}{n17}{\tiny 1}
\edge{n17}{n18} \edgetext{n17}{n18}{\tiny 1}
\edge{n18}{n19} \edgetext{n18}{n19}{\tiny 1}
\edge{n19}{n20} \edgetext{n19}{n20}{\tiny 1}
\edge{n20}{n21} \edgetext{n20}{n21}{\tiny 1}
\edge{n22}{n23} \edgetext{n22}{n23}{\tiny 1}
\edge{n23}{n24} \edgetext{n23}{n24}{\tiny 1}
\edge{n24}{n25} \edgetext{n24}{n25}{\tiny 1}
\edge{n25}{n26} \edgetext{n25}{n26}{\tiny 1}
\edge{n26}{n27} \edgetext{n26}{n27}{\tiny 1}
\edge{n27}{n28} \edgetext{n27}{n28}{\tiny 1}
\edge{n29}{n30} \edgetext{n29}{n30}{\tiny 1}
\edge{n30}{n31} \edgetext{n30}{n31}{\tiny 1}
\edge{n31}{n32} \edgetext{n31}{n32}{\tiny 1}
\edge{n32}{n33} \edgetext{n32}{n33}{\tiny 1}
\edge{n33}{n34} \edgetext{n33}{n34}{\tiny 1}
\edge{n34}{n35} \edgetext{n34}{n35}{\tiny 1}
\edge{n36}{n37} \edgetext{n36}{n37}{\tiny 2}
\edge{n37}{n38} \edgetext{n37}{n38}{\tiny 1}
\edge{n38}{n39} \edgetext{n38}{n39}{\tiny 2}
\edge{n39}{n40} \edgetext{n39}{n40}{\tiny 2}
\edge{n40}{n41} \edgetext{n40}{n41}{\tiny 2}
\edge{n41}{n42} \edgetext{n41}{n42}{\tiny 2}

\edge{n1}{n8} \edgetext{n1}{n8}{\tiny 3}
\edge{n2}{n9} \edgetext{n2}{n9}{\tiny 2}
\edge{n3}{n10} \edgetext{n3}{n10}{\tiny 2}
\edge{n4}{n11} \edgetext{n4}{n11}{\tiny 1}
\edge{n5}{n12} \edgetext{n5}{n12}{\tiny 1}
\edge{n6}{n13} \edgetext{n6}{n13}{\tiny 1}
\edge{n7}{n14} \edgetext{n7}{n14}{\tiny 3}
\edge{n8}{n15} \edgetext{n8}{n15}{\tiny 1}
\edge{n9}{n16} \edgetext{n9}{n16}{\tiny 1}
\edge{n10}{n17} \edgetext{n10}{n17}{\tiny 1}
\edge{n11}{n18} \edgetext{n11}{n18}{\tiny 2}
\edge{n12}{n19} \edgetext{n12}{n19}{\tiny 2}
\edge{n13}{n20} \edgetext{n13}{n20}{\tiny 2}
\edge{n14}{n21} \edgetext{n14}{n21}{\tiny 1}
\edge{n15}{n22} \edgetext{n15}{n22}{\tiny 3}
\edge{n16}{n23} \edgetext{n16}{n23}{\tiny 2}
\edge{n17}{n24} \edgetext{n17}{n24}{\tiny 2}
\edge{n18}{n25} \edgetext{n18}{n25}{\tiny 1}
\edge{n19}{n26} \edgetext{n19}{n26}{\tiny 1}
\edge{n20}{n27} \edgetext{n20}{n27}{\tiny 1}
\edge{n21}{n28} \edgetext{n21}{n28}{\tiny 3}
\edge{n22}{n29} \edgetext{n22}{n29}{\tiny 1}
\edge{n23}{n30} \edgetext{n23}{n30}{\tiny 1}
\edge{n24}{n31} \edgetext{n24}{n31}{\tiny 1}
\edge{n25}{n32} \edgetext{n25}{n32}{\tiny 2}
\edge{n26}{n33} \edgetext{n26}{n33}{\tiny 2}
\edge{n27}{n34} \edgetext{n27}{n34}{\tiny 2}
\edge{n28}{n35} \edgetext{n28}{n35}{\tiny 1}
\edge{n29}{n36} \edgetext{n29}{n36}{\tiny 3}
\edge{n30}{n37} \edgetext{n30}{n37}{\tiny 2}
\edge{n31}{n38} \edgetext{n31}{n38}{\tiny 2}
\edge{n32}{n39} \edgetext{n32}{n39}{\tiny 1}
\edge{n33}{n40} \edgetext{n33}{n40}{\tiny 1}
\edge{n34}{n41} \edgetext{n34}{n41}{\tiny 1}
\edge{n35}{n42} \edgetext{n35}{n42}{\tiny 3}

\end{graph}
\caption{An interior point of $5\P$ for $\G(6,7)$.}
\label{mevennodd5interior}
\end{figure}
Since this graph labelling is not symmetric, we can get another
interior point by ``flipping'' the graph (i.e., the weight for edges
$\{(i,j), (i,j+1)\}$ and $\{(i,j), (i+1,j)\}$ are the weights for
$\{(n-1-i,j), (n-1-i,j+1)\}$ and $\{(n-2-i,j), (n-1-i,j)\}$,
respectively, from the interior point given previously). Thus, $5\P$
does not have a unique interior point so $\P$ is not Gorenstein.
\end{proof}

We give an example of the Ehrhart polynomial and Ehrhart series of the
perfect matching polytope of a grid graph. We used polymake
\cite{polymake} for conversion from the hyperplane description given
in Theorem \ref{hyperplane} to vertex descriptions, and normaliz
\cite{normaliz} for the computation of the Ehrhart functions.

\begin{example} \rm The Ehrhart functions of $\P(3,4)$ are
  \begin{align*}
    L_{\P(3,4)}(t) &= \frac{1}{120} (t+1)(t+2)(t+3)(t+4)(t^2 + 5t +
    5)\,,\\
    \Ehr_{\P(3,4)}(z) &= \frac{z^2 + 4z + 1}{(1-z)^7}\,.
  \end{align*}
  The smallest example of a non-Gorenstein polytope arises from
  $\P(4,5)$, which has Ehrhart series
  \[ 
  \Ehr_{\P(4,5)}(z) = \frac{21z^8 + 760z^7 + 5919z^6 + 15578z^5 +
    16432z^4 + 7356z^3 + 1339z^2 + 82z + 1}{(1-z)^{13}}\,. 
  \]
\end{example}


\section{Unimodular Triangulations and
  Unimodality}\label{unimodalsection}

We remind the reader that a {\bf triangulation} of a polytope $\P$ is
a set of simplices $\Delta$ whose union is $\P$ such that the
intersection of any two simplices is a face of both. An integral
simplex with vertices $v_0, \dots, v_n$ is {\bf unimodular} if
$\{v_1-v_0, \dots, v_n-v_0\}$ forms a basis of $\Z^n$ as a free
Abelian group, and a triangulation is unimodular if each simplex is
unimodular.

One might ask if $\P(m,n)$ has a unimodular
triangulation. Fortunately, the answer is yes. Given an ordering $\tau
= (v_1, \dots, v_n)$ of the vertices of a polytope $\P$, the {\bf
  reverse lexicographic triangulation} $\Delta_\tau$ is defined as
follows. For a single vertex $v$, $\Delta_\tau = \{v\}$.  In general,
consider all facets (codimension 1 faces) of $\P$ that do not contain
$v_n$. For each such facet $F$, there is an ordering of the vertices
of $F$ induced by $\tau$ that gives a reverse lexicographic
triangulation of $F$. For each simplex in this triangulation, take the
convex hull of its union with $v_n$. The union of all such convex
hulls is $\Delta_\tau$. A polytope is said to be {\bf compressed} if
the reverse lexicographic triangulation with respect to any ordering
of its vertices is unimodular. The following result immediately
guarantees the existence of unimodular triangulations.

\begin{theorem}[Santos~\cite{paco}, Ohsugi-Hibi~\cite{revlex},
  Sullivant~\cite{SullivantCompressed}]
  \label{ohsugihibi}
  If $\P$ is a 0/1 polytope
  (i.e., each coordinate of each vertex is either 0 or 1) defined by
  the linear equalities and inequalities
  \[ 
  {\bf A \, x} = {\bf u} \quad \text{and} \quad 0 \le x_i \le 1\,,
  \] 
  where ${\bf A}$ is an integer-valued matrix and ${\bf u}$ is an
  integer-valued column vector, then $\P$ is compressed.
\end{theorem}

As a bonus, we obtain Theorem~\ref{unimodalthm}, that is, for the
$(m,n)$ discussed in Theorem \ref{gorensteinthm}, the Ehrhart
$h$-vector is unimodal. This is an immediate consequence of
Theorem~\ref{ohsugihibi} combined with the following result.

\begin{theorem}[Athanasiadis--Bruns--R\"omer
  \cite{compressed,gorensteinhvector}]
  \label{compressedthm} If
  $\P$ is compressed and Gorenstein, then the Ehrhart $h$-vector of
  $\P$ is unimodal.
\end{theorem}


\section{Domino Stackings and Magic Labellings} \label{dominosection}

There is a natural many-to-one function from domino stackings of
height $t$ of an $m \times n$ rectangular board to the magic
labellings of sum $t$ of $\G(m,n)$. It is conceivable that some magic
labelling cannot be realized as a physical stacking. However, we have
shown earlier that $\P(m,n)$ has a unimodular triangulation, which
implies that it is a normal polytope. (An integral polytope $\P$ is
{\bf normal} if every integer point in $t\P$ can be written as a sum
of $t$ integer points in $\P$.) In particular, this means that every
magic labelling of sum $t$ of $\G(m,n)$ must be the sum of $t$ perfect
matchings, so our function is a surjection, and this is what we
claimed in Proposition~\ref{normalprop}.

Recall that for fixed $m$, the generating function $\sum_{n \ge 0}
T(m,n,1) z^n$ evaluates to a proper rational function in $z$, which
means that the sequence $(T(m,n,1))_{n \ge 0}$ is given by a linear
recurrence relation. Running these recurrence relations backwards, one
obtains values for $T(m,n,1)$ when $n<0$. It was shown in
\cite{dominoreciprocity} that these values for $n<0$ are the same even
if the recurrence relation used is not minimal. A reciprocity relation
was also given, which is the content of Theorem \ref{recurrthm}.

This is equivalent to a result in \cite{dimercovering} which gives
properties of the numerator and denominator 
polynomials in the rational function
$\sum_{n \ge 0} T(m,n,1) z^t$. We would like an analogous result for
magic labellings of sum $t$ and domino stackings of height $t$. Recall
that Theorem~\ref{recurrence} states that for general $t$, if one
fixes $m$, then there is a recurrence relation in $n$ for the function
$T(m,n,t)$. The proof that we give is essentially the same as the one
found in \cite{dominotilings} with only 
slight modifications.

\begin{proof}[Proof of Theorem~\ref{recurrence}]
  Given a magic labelling of sum $t$ for $\G(m,n)$, consider its
  canonical planar drawing. We can uniquely encode it with $n+1$
  vectors $c_0, \dots, c_n$ of length $m$ in the following way. For $1
  \le i \le n-1$, the $j\th$ entry of $c_i$ is the label of the $j\th$
  horizontal edge counting from the top of the graph in the $i\th$
  column starting from the left and going right. We define $c_0$ and
  $c_n$ to be zero vectors of length $m$.

  We now form a directed graph $\mathscr{G}_m$ whose vertex set is
  $\{0, 1, \dots, t\}^m$. For two vectors $u,v \in \{0, \dots, t\}^m$,
  we form the edge $(u,v)$ if and only if there is some $n$ such that
  there is a magic labelling of sum $t$ for $\G(m,n)$ encoded by the
  vectors $c_0, \dots, c_n$ with $c_i = u$ and $c_{i+1} = v$ for some
  $i$. Given this, there is a natural bijection between the magic
  labellings of sum $t$ of $\G(m,n)$ and walks of length $n$ in
  $\mathscr{G}_m$ from $\{0\}^m$ to itself.

  Order the vertices of $\mathscr{G}_m$ by interpreting each vector as
  a number in base $t+1$, and let $A_m$ be the adjacency matrix of
  $\mathscr{G}_m$. The number of walks of length $n$ from $\{0\}^m$ to
  itself is the $(0,0)$ entry of $A_m^n$. By the Cayley--Hamilton
  theorem, there is a polynomial relation
  \[ 
  A_m^d + b_{d-1} A_m^{d-1} + \cdots + b_1 A_m + b_0 = 0 
  \] 
  for $d = (t+1)^m$ and integers $b_i$. Multiplying both sides by
  $A_m^{n-d}$ and looking at the $(0,0)$ entry in each term gives a
  recurrence relation for $T(m,n,t)$ involving the coefficients $b_i$.

  In the case of domino stackings of height $t$, we can get a similar
  result. Instead of vectors in $\{0,\dots,t\}^k$ making up the
  vertices of $\mathscr{G}_m$, we encode each domino stacking as a $k
  \times t$ zero-one matrix and play the same game to get a recurrence
  relation for the sequence $(T(m,n,1)^t)_{n \ge 0}$.
\end{proof}

As an immediate corollary, we now know that for fixed $m$ and $t$, the
generating function $\sum_{n \ge 0} T(m,n,t) z^n$ evaluates to a
rational function $\frac{P(z)}{Q(z)}$. The next step is to understand
what properties $P$ and $Q$ have.

The sequence $(T(2,n,1))_{n \ge 0}$ is a shift of the Fibonacci
numbers, so this gives a nice proof that powers (for a fixed exponent)
of Fibonacci numbers can be encoded by a linear homogeneous recurrence
relation. This was previously done in \cite{fibonaccipowers}, which
also gives recurrence relations for the {\it coefficients} of the
recurrence relations.

\section{Convex Geometric Considerations of Domino Tiling
  Reciprocity}\label{geometricsection}

In this section, we give a convex geometric proof of why reciprocity
exists for the number of domino tilings of the $m \times n$
rectangular board for fixed $m$ using the following theorem which
gives a closed form. We make some remarks on reciprocity of domino
stackings of fixed height $t$ and width $m$ at the end of the proof.

\begin{theorem}[Kasteleyn \cite{closedform}] \label{closedformtheorem}
  The number of domino tilings of the $m \times n$ rectangular board
  is given by
  \[
  T(m,n,1) = \prod_{j=1}^{\lceil m/2 \rceil} \frac{c_j^{n+1} -
    \bar{c}_j^{n+1}}{c_j - \bar{c}_j}
  \]
  where
  \begin{align*}
    c_j &= \cos \frac{j\pi}{m+1} + \sqrt{1 + \cos^2
      \frac{j\pi}{m+1}}\,,\\
    \bar{c}_j &= \cos \frac{j\pi}{m+1} - \sqrt{1 + \cos^2
      \frac{j\pi}{m+1}}\,.
  \end{align*}
\end{theorem}

We remark here that the above theorem is in a different form from the
one found in \cite{closedform}. The original formulation gives a
product from $j=1$ to $j=\lfloor m/2 \rfloor$, and designates
$T(m,n,1) = 0$ when $mn$ is odd. However, when $mn$ is odd, $c_j^{n+1}
- \bar{c}_j^{n+1} = 0$ for $j = \lceil m/2 \rceil$, and for $m$ odd
and $n$ even, it is $c_j - \bar{c}_j$, so our formulation is
equivalent.

First note that $c_j\overline{c}_j = -1$ for all $j$. The key
observation is that
\[
\frac{c_j^{n+1} - \bar{c}_j^{n+1}}{c_j - \bar{c}_j} = c_j^n +
c_j^{n-1}\bar{c}_j + \cdots + c_j\bar{c}_j^{n-1} + \bar{c}_j^n\,,
\]
and that the sum on the right is an evaluation of a certain
multivariate generating function \eqref{eq:gf} below.  That is,
consider the polytope $\P = \conv\{(1,0), (0,1)\}$. Let $\ell :=
\lceil m/2 \rceil$; we are interested in the $\ell$-fold Cartesian
product of $\P$ with itself. Let $\CC \subseteq \R^{2\ell+1}$ be the
cone over $\P$. To be precise, $\CC = \{ (x,\lambda) \in \R^\ell
\times \R : \lambda \ge 0,\ x \in \lambda P \}$.  We consider the {\bf
  integer point transform} of $\CC$, which is the following formal
Laurent series in $2\ell+1$ variables.
\[ \tag{$*$} \label{eq:gf} \sigma_{\CC}({\bf x}) := \sum_{{\bf y} \in
  \CC \cap \Z^{2\ell+1}} {\bf x}^{{\bf y}}\,,
\]
where ${\bf x}^{{\bf y}} := x_1^{ y_1 } x_2^{ y_2 } \cdots x_{ 2 \ell
  + 1 }^{ y_{ 2 \ell + 1 } }$. In the case of a simplicial cone
$\CC'$, if $y_1, \dots, y_k$ are the generators of $\CC'$, define the
{\bf fundamental parallelepiped}
\[
\Pi(\CC') := \{ \lambda_1y_1 + \cdots + \lambda_ky_k : 0 \le \lambda_i
< 1 \}\,.
\]
Then the integer point transform $\sigma_{\CC'}({\bf x})$ is a
rational function of the form (see, e.g., \cite[Chapter 3]{ccd})
\[
\sigma_{\CC'}({\bf x}) = \frac{\sigma_{\Pi(\CC')}({\bf
    x})}{(1-{\bf x}^{y_1}) \cdots (1-{\bf x}^{y_k})}\,.
\]
We can triangulate $\P$ (using no new vertices), and
this triangulation induces a triangulation of $\CC$ into simplicial
cones. The generators of each will be of the form $x_{i_1}x_{i_2}
\cdots x_{i_{\ell}}$ where $\{i_1, i_2, \dots, i_{\ell}\} \subseteq
\{1, \dots, 2\ell\}$
contains exactly one of $2i-1$ and $2i$ for each
$i=1, \dots, \ell$.
Conversely, each such subset will appear as the
generator for at least one of the simplicial cones in the
triangulation of $\CC$. By adding the rational functions for these
simplicial cones and playing inclusion-exclusion on the intersections,
the integer point transform of $\CC$ is
\[ 
\sigma_{\CC}({\bf x}) = 
\frac{P({\bf x})}{\prod (1-x_{i_1}x_{i_2} \cdots x_{i_{\ell}})}
\]
where the product runs over all $\ell$-subsets of $\{1, \dots,
2\ell\}$ that contain exactly one of $2i-1$ and $2i$ for each $i=1,
\dots, \ell$, and $P$ is some multivariate polynomial. We make the
substitutions $x_{2i-1} \mapsto c_i$, $x_{2i} \mapsto \bar{c}_i$, and
$x_{2\ell+1} \mapsto z$. The polytope $n\P$ can be described as the
set of nonnegative solutions to the equations $x_1 + \cdots +
x_{2\ell} = \ell n$ and $x_{2i-1} + x_{2i} = n$ for each $i=1, \dots,
\ell$. Thus, using Theorem~\ref{closedformtheorem},
\begin{align*}
  \sigma_{n\P}(c_1, \bar{c}_1, \dots, c_{\ell}, \bar{c}_{\ell}, z) &=
  z^n \prod_{j=1}^{\ell} (c_j^n + c_j^{n-1}\bar{c}_j + \cdots +
  c_j\bar{c}_j^{n-1} + \bar{c}_j^n)\\
  &= T(m,n,1)z^n\,.
\end{align*}
From this substitution,
\begin{align*}
  \sigma_{\CC}(c_1, \bar{c}_1, \dots, c_{\ell}, \bar{c}_{\ell}, z) &=
  \sum_{n \ge 0} T(m,n,1) z^n\\
  &= \frac{P(c_1, \bar{c}_1, \dots, c_{\ell}, \bar{c}_{\ell},
    z)}{\prod_I (1-d_1^{i_1} d_2^{i_2} \cdots
    d_{\ell}^{i_{\ell}}z)}\,,
\end{align*}
where $I$ runs over all possible binary strings of length $\ell$, and
$i_j$ denotes the $j\th$ bit of $I$, and $d_j^0 = c_j$ and $d_j^1 =
\bar{c}_j$.

To proceed, we use the following theorem due to Stanley which gives a
relation between the integer point transform of a rational cone and
the integer point transform of its interior.

\begin{theorem}[Stanley \cite{stanleyreciprocity}] If $\CC$ is a
  rational cone in $\R^d$ whose apex is the origin, then as rational functions,
  \[
  \sigma_{\CC}(z_1^{-1}, \dots, z_d^{-1}) = (-1)^{\dim \CC}
  \sigma_{\CC^{\circ}} (z_1, \dots, z_d)\,,
  \]
  where $\CC^{\circ}$ denotes the relative interior of $\CC$.
\end{theorem}

Equipped with this fact, we have
\begin{align} \begin{split} \label{reciprocitystep}
    \sigma_{\CC^{\circ}}(c_1, \bar{c}_1, \dots, c_{\ell},
    \bar{c}_{\ell}, z) &= (-1)^{\ell+1} \sigma_{\CC}(-\bar{c}_1, -c_1,
    \dots,
    -\bar{c}_{\ell}, -c_{\ell}, z^{-1})\\
    &= (-1)^{\ell+1} \sigma_{\CC}(c_1, \bar{c}_1, \dots, c_{\ell},
    \bar{c}_{\ell}, (-1)^{\ell} z^{-1})\\
    &= (-1)^{\ell+1} \sum_{n \le 0} (-1)^{\ell n} T(m,-n,1) z^n\\
    &= (-1)^{\ell} \sum_{n > 0} (-1)^{\ell n} T(m,-n,1) z^n .
\end{split} \end{align}
The last equality follows from the identity (of rational functions)
\[
\sum_{n \ge 0} f(n) z^n + \sum_{n < 0} f(n) z^n = 0\,,
\]
where $f(n)$ is given by a linear recurrence relation. Finally, we use
the fact that $\CC$ is a Gorenstein cone 
(i.e., adding some
 vector to the cone allows us to enumerate the
integer points in its interior). This means that 
\begin{align} \begin{split} \label{interiorstep}
  \sigma_{\CC^{\circ}}(c_1, \bar{c}_1, \dots, c_{\ell},
  \bar{c}_{\ell}, z) &= c_1 \bar{c}_1 \cdots c_{\ell} \bar{c}_{\ell}
  z^2 \sigma_{\CC}(c_1, \bar{c}_1, \dots,
  c_{\ell}, \bar{c}_{\ell}, z)\\
  &= (-1)^{\ell} z^2 \sum_{n \ge 0} T(m,n,1) z^n .
\end{split} \end{align}
Putting both sets of equalities together and matching coefficients
gives 
\[
T(m,n,1) = (-1)^{\ell n} \, T(m,-n-2,1)\,.
\]
To recover the result in \cite{dominoreciprocity}, note that when $m
\equiv 3 \pmod 4$ or $m \equiv 0 \pmod 4$, $\ell$ is even, so
\[
T(m,n,1) = T(m,-n-2,1)\,.
\]
In the other cases, $\ell$ is odd, so 
\[
T(m,n,1) = (-1)^n \, T(m,-n-2,1)\,.
\] 
In particular, for $m \equiv 1 \pmod 4$, $T(m,n,1) = 0$ when $n$ is
odd, so this is equivalent to 
\[
T(m,n,1) = T(m,-n-2,1)\,, \quad m \equiv 1 \pmod 4\,,
\]
which gives the desired reciprocity relations.

For general $t$, raise both sides of the reciprocity equations to the
$t\th$ power to get reciprocity relations for domino stackings of
height $t$. Alternatively, we could run through the above proof with
the following small modifications. Consider the $t$-fold product of
$\P$ with itself. The cone over $\P^t$ is now $(t\ell+1)$-dimensional,
so the last term in equation (\ref{reciprocitystep}) becomes
\[
(-1)^{t\ell} \sum_{n > 0} (-1)^{t\ell n} \, T(m,-n,1)^t z^n ,
\]
and the last term in equation (\ref{interiorstep}) becomes
\[
(-1)^{t\ell} z^2 \sum_{n \ge 0} T(m,n,1)^t z^n .
\]
Combining these two gives the general reciprocity relation
\[
T(m,n,1)^t = (-1)^{t\ell n} \, T(m,-n-2,1)^t .
\]


\section{Further Directions}

A natural question to ask is if unimodality of the Ehrhart $h$-vector
extends to the $(m,n)$ such that $\P(m,n)$ is not
Gorenstein. Unfortunately, computing examples larger than $\G(4,5)$
becomes difficult, so there is not much evidence supporting whether
they are or are not unimodal.

We defined $\G(m,n)$ having vertex set $\{(i,j) \in \Z^2 : 0 \le i <
m\,,\ 0 \le j < n\}$ with $\{(i,j), (i',j')\}$ an edge if and only if
$|i-i'| + |j-j'| = 1$, but there is no reason to restrict to
2-dimensional lattices. We can define $k$-dimensional grid graphs
$\G(m_1,m_2, \dots, m_k)$ having vertex set $\{(x_1, \dots, x_k) \in
\Z^k : 0 \le x_i < m_i \}$ such that $\{(x_i), (y_i)\}$ is an edge if
and only if $\sum_{i=1}^k |x_i - y_i| = 1$. What further things can
one say about these higher dimensional analogues? It would make sense
that very few of these graphs have Gorenstein polytopes.  The perfect
matching polytope of the cube graph
$\G(2,2, \dots, 2)$ is Gorenstein for arbitrary
$k$, since this is a $k$-regular graph. Some questions one can ask is
if there is a largest $k$ such that $\G(2,2, \dots, 2)$ is the only
graph with a Gorenstein polytope (assuming that $m_i > 1$). If not, is
there a largest $k$ such that the number of nontrivial Gorenstein
polytopes $\P(m_1, \dots, m_k)$ is infinite? We can also ask about the
function $T(m_1, \dots, m_k; t)$ that counts the number of magic
labellings of sum $t$ of $\G(m_1, \dots, m_k)$. Since
multi-dimensional grid graphs are bipartite, $\sum_{t \ge 0} T(m_1,
\dots, m_k; t) z^t$ evaluates to a rational function. It is not
difficult to see using the proof of Theorem \ref{recurrence} that
$\sum_{m_i \ge 0} T(m_1, \dots, m_k; t) z^{m_i}$ is also a rational
function.

One could also define two domino stackings to be equivalent if the
layers of tilings of one is a permutation of the layers of tilings of
the other. Is there a recurrence relation for this new set of
stackings? New techniques would be needed as the adjacency matrix
argument no longer works.

As mentioned earlier, Riordan gives recurrence relations for
coefficients of recurrence relations for the sequences given by powers
of Fibonacci numbers in \cite{fibonaccipowers}. It is conceivable that
the same phenomenon occurs for powers of $T(m,n,1)$ for general
$m$. Different techniques from what is presented here would most
likely have to be employed. The main problem arises in that the proof
of Theorem \ref{recurrence} does not give minimal recurrence
relations. For instance, in the case $m=2$ and $t=1$, which is a shift
of the Fibonacci sequence, a 5-term recurrence relation is given by
the proof, but only a 3-term recurrence relation is needed. This
corresponds to the fact that if this non-minimal recurrence relation
appears in the denominator of the rational function $\sum_{n \ge 0}
T(m,n,1)^t z^n$, then it will not be in reduced terms. A further
problem is that the adjacency matrix argument given above presents a
recurrence relation that works for {\it all} of the entries of the
matrix, but we are only interested in one specific entry.

Along those same lines, for fixed $t$, the sequence $T(m,n,t)$ for
counting magic labellings of $\G(m,n)$ is given by a linear
homogeneous recurrence relation if either $m$ or $n$ is fixed. It
would be interesting to find a recurrence relation when $m$ and $n$
are both allowed to vary, and also to see if the coefficients that
arise from these recurrence relations can also be shown to satisfy
their own recurrence relation.


\bibliographystyle{amsplain}

\def\cprime{$'$} \def\cprime{$'$}
\providecommand{\bysame}{\leavevmode\hbox to3em{\hrulefill}\thinspace}
\providecommand{\MR}{\relax\ifhmode\unskip\space\fi MR }
\providecommand{\MRhref}[2]{%
  \href{http://www.ams.org/mathscinet-getitem?mr=#1}{#2}
}
\providecommand{\href}[2]{#2}

\end{document}